\documentclass{amsart}
\usepackage{graphicx}
\usepackage{latexsym}
\usepackage{hyperref}

\usepackage{amsfonts}
\usepackage[all]{xy}

\newtheorem{corollary}{Corollary}

\newtheorem{lemma}{Lemma}[section]

\newtheorem{theorem}{Theorem}

\numberwithin{equation}{section}

\theoremstyle{definition}

\numberwithin{equation}{section} \theoremstyle{remark}

\newcommand{\be}{\begin{eqnarray}}
\newcommand{\en}{\end{eqnarray}}
\newcommand{\no}{\nonumber}

\newcommand{\de}{\delta}

\newcommand{\fr}{\frac}
\newcommand{\pa}{\partial}

\newcommand{\D}{\Delta}
\newcommand{\om}{\Omega}
\newcommand{\na}{\nabla}
\newcommand{\R}{\mathbb{R}}
\newcommand{\lan}{\langle}
\newcommand{\ra}{\rangle}

\newcommand{\ri}{\rightarrow}

\newcommand{\spacing}[1]{\renewcommand{\baselinestretch}{#1}\large\normalsize}
\spacing{1.3}
\begin{document}
\title [On Payne-Schaefer's Conjecture about an Overdetermined Boundary Problem of Sixth
order]{On Payne-Schaefer's Conjecture about an Overdetermined
Boundary Problem of Sixth order }
\author[Changyu Xia] {Changyu Xia}

\maketitle \thispagestyle{empty}
\begin{abstract}
This paper considers overdetermined boundary problems. Firstly, we
give a proof to the Payne-Schaefer conjecture about an
overdetermined  problem of sixth order in
 the two  dimensional case and under an additional condition for the
 case of dimension no less than three. Secondly, we  prove an integral
 identity for an overdetermined problem of fourth order which can be
 used to deduce Bennett's symmetry theorem. Finally, we prove a
 symmetry result for an overdetermined problem of second order by
 integral identities.
\end{abstract}


\section{Introduction and the main results}

\markboth{\hfill
  Changyu Xia\hfill}
{\hfill On Payne-Schaefer's Conjecture about an Overdetermined
Boundary Problem \hfill}

In a celebrated paper in 1971, Serrin  initiated the study of
elliptic equations under overdetermined boundary condition and
established in particular the following  seminal result.
\begin{theorem}\label{th1}(\cite{s}). If $\om$ is a bounded domain with smooth
boundary in $\mathbb{R}^n$ and if the solution to the problem \be
\left\{\begin{array}{l}
\Delta u = -1 \ \ \ {\rm in  \ }\ \ \om,\\
u =  0 \ \ \ \ \ \ \ \ {\rm on  \ }\ \pa\om,
\end{array}\right.
\en has the property that $\pa u/\pa\nu$ is equal to a constant $c$
on $\pa\om$, then $\om$ is a ball of radius $|nc|$ and
$u=(n^2c^2-r^2)/2n$, where $\nu$ is the outward unit normal of
$\pa\om$ and  $r$ is the distance from the center of the ball.
\end{theorem}
Several proofs to the above result have appeared. Serrin's proof is
based on the Hopf maximum principle  and a
reflection-in-moving-planes argument which could be extended to more
general elliptic equations and somewhat more general boundary
conditions. A simple proof of Serrin's result based on a Rellich
identity and a maximum principle was given by Weinberger \cite{w}.
By the method of duality theorem Payne and Schaefer \cite{ps} gave a
proof of Theorem 1 which does not make explicit use of maximum
principle. Choulli and Henrot \cite{ch} used domain derivative to
prove Serrin's theorem which also does not use the maximum principle
explicitly. From the need to extend Serrin overdetermined result to
non uniformly elliptic operators of Hessian type, Brandolini,
Nitsch, Salani and Trombetti \cite{bnst1} used an integral approach
via arithmetic-geometric mean inequality to prove Serrin's theorem
and they also established the stability of the Serrin problem
\cite{bnst2}. Serrin's theorem is a landmark in the study of
overdetermined boundary value problem. The ideas and techniques in
proving Serrin's theorem have been widely used and generalized to
prove symmetry for more general  overdetermined problems. Troy
\cite{t} used Serrin's moving planes method to prove a symmetry
theorem for a system of semilinear elliptic equations, Alessandrini
\cite{a} adapted this method to  condensers in a capacity problem.
Garofalo and Lewis \cite{gl} extended Weinberger's method to more
general second order partial differential equations. In \cite{ps},
Payne and Schaefer  studied overdetermined problems of higher
orders, obtained various symmetry results and proposed the following
important \vskip0.3cm \noindent{\bf Conjecture (\cite{ps}).} {\it
Let $\om$ be a bounded domain in $\R^n$ with smooth boundary, If u
is a sufficiently smooth solution of the following overdetermined
problem: \be\label{c}
 \left\{\begin{array}{l}
  \D^3 u=-1 \ \ \ \ \ \ \ \ \ \ \ \ \hskip0.1cm {\rm in}\ \ \  \om, \\
  u=\fr{\pa u}{\pa \nu}=\D u=0\ \ \ {\ \rm  on\ } \pa\om, \\
  \fr{\pa (\D u)}{\pa \nu}= c \ \ \ \ \ \ \ \ \ \ \ \ \hskip0.05cm {\ \rm on\ }
  \pa\om,
\end{array}\right.
\en then $\om$ is an $n$-ball.} \vskip0.3cm In this paper, we prove
Payne-Schaefer's conjecture in the case $n=2$ and also prove the
case $n\geq 3$ under an additional hypothesis.

\begin{theorem}\label{th1} Let $\om$ is a bounded domain in $\R^n, n\geq 2$, with $C^{6+\epsilon}$ boundary
and if the following overdetermined problem has a solution in
$C^6(\overline{\om})$: \be
& & \D^3 u=-1 \ \ \ \ \ \ \ \ \ \ \ \ \ {\rm in}\ \ \  \om, \\
& &   u=\fr{\pa u}{\pa \nu}=\D u=0 \ \ {\ \rm on\ } \pa\om,  \\
& &   \fr{\pa (\D u)}{\pa \nu}= c \ \ \ \ \ \ \ \ \ \ \hskip0.1cm {\
\rm on\ } \pa\om, \en  where  $c$ is a constant. When $n\geq 3$,
assume that \be\label{th1.1} \int_{\om} (\D^2u)^2\gamma dx\leq
\fr{2(n+2)c^2|\om|}{n+6}.\en Here $|\om|$ denotes the volume of
$\om$ and $\gamma$ is the torsion function of $\om$ given by
\be\no\left\{\begin{array}{l} \D \gamma = -1 \ \ {\rm in\ \ \ } \om,\\
\gamma = 0 \ \ \ \ \ \  {\rm on} \ \ \pa\om.
\end{array} \right. \en Then $\om$ is a ball of radius $\left(|c|n(n+2)(n+4)\right)^{\fr 13}$, and
\be\no u(x) &=&  -\fr 1{48n(n+2)(n+4)}r^6
+\left(\fr{c^2}{n(n+2)(n+4)}\right)^{\fr 13}\cdot\fr{r^4 }{16}\\
 & & -\left(c^4n(n+2)(n+4)\right)^{\fr
13}\cdot\fr{r^2}{16}+\fr{c^2n(n+2)(n+4)}{48},\en where $r$ denotes
the distance from $x$ to the  center of $\om$.
\end{theorem}
It should be mentioned that for a ball in $\R^n$, (\ref{th1.1})
becomes an equality. We shall explain this in the next section.

An integral dual for (1.3)-(1.5) is \be \int_{\om} \phi dx=
c\int_{\pa\om} \D\phi ds \en for any triharmonic function $\phi$ in
$\om$ for which $\phi=\pa\phi/\pa\nu=0$ on $\pa\om$. Thus, we have
from Theorem 2 the following \vskip0.3cm\begin{corollary}
 Let $\om$ is a bounded domain in $\R^2$,  with $C^{6+\epsilon}$ boundary
and if (1.8) holds for any triharmonic function $\phi$ in $\om$ for
which $\phi=\pa\phi/\pa\nu=0$ on $\pa\om$, where $c$ is a constant.
Then $\om$ is a disk.
\end{corollary}

In \cite{b}, Bennett established the following  symmetry result.

\begin{theorem}
If $\om$ is a bounded domain in $\R^n$ with $C^{4+\epsilon}$
boundary and if the following overdetermined problem has a solution
in $C^4(\overline{\om})$: \be \left\{\begin{array}{l}
\Delta^2 u = -1 \ \ \ \ \ \ {\rm in  \ }\ \ \om,\\
u = \fr{\pa u}{\pa\nu}= 0 \  \ \hskip0.04cm \ \ {\rm on  \ }\ \pa\om,\\
\D u\equiv c \ \ \ \ \ \ \ \hskip0.07cm \ \ {\rm on \ } \ \pa\om\
(c\ constant),
\end{array}\right.
\en then $\om$ is a ball of radius $(|c|n(n+2))^{1/2}$, and \be
u(x)=\fr{-1}{2n}\left\{\fr 14(n+2)(nc)^2+\fr{nc}2 r^2+\fr
1{4(n+2)}r^4\right\}, \en where where $r$ denotes the distance from
x to the  center of $\om$. \end{theorem}

A crucial point in  Bennett's proof is to use the following identity
\cite{p}: \be\no
 \fr 12\D\Phi &=&  \sum_{i,j,k}u_{ijk}^2-\fr 3{n+2}|\na(\D u)|^2\\
&=& \sum_{i,j,k}\left\{u_{ijk}-\fr 1{n+2}\left((\D
u)_i\delta_{jk}+(\D u)_j\delta_{ik}+(\D
u)_k\delta_{ij}\right)\right\}^2, \en where \be
\Phi=\fr{n-4}{n+2}u+\fr{n-4}{2(n+2)}(\D u)^2+|\na^2u|^2-\lan\na u,
\na(\D u)\ra. \en Since \be \Phi|_{\pa\om}=\fr{3nc^2}{2(n+2)},\en it
then follows from the maximum principle that \be \Phi\leq
\fr{3nc^2}{2(n+2)},\ \ \ {\rm in} \ \om.\en On the other hand, from
Green's theorem and Rellich identity, one has \be \int_{\om} \Phi
dx=\fr{3nc^2}{2(n+2)}\cdot|\om|. \en Thus, $\Phi\equiv
3nc^2/(2(n+2))$ in $\overline{\om}$ and so $\D\Phi\equiv 0$ in
$\overline{\om}$. Therefore, each term of the sum on the right hand
side of (1.11) vanishes which implies  that \be (\D u)_{ij}=-\fr
1n\delta_{ij}. \en One can then obtain  the conclusions of Theorem 3
easily.

 In  this paper, we  obtain an
 integral identity for an overdtermined problem of fourth order from which
 one can prove  Bennett's theorem without using the subharmonicity of the function $\Phi$.
\begin{theorem} Let $\om$ is a bounded domain in $\R^n, n\geq 2 $ with
$C^{4+\epsilon}$ boundary. Let $g: \R\ri\R$ be a $C^2$ function and
set $G(t)=\int_0^t g(s)ds$. If $u\in C^4(\overline{\om})$ is a
solution of the following overdetermined problem : \be & & \D^2
u=-g(u) \ \ \ \ \ \  {\rm in}\ \ \
 \om,\\
& &   u=\fr{\pa u}{\pa \nu}=0 \ \ \ \ \ \ \ {\ \rm on\ } \pa\om, \\
& & \D u = c \ \ \ \ \ \ \ \ \ \ \ \ \hskip0.1cm {\ \rm on\ }
\pa\om, \en where  $c$ is a constant. Then, we have \be\label{th2.0}
\no & & \int_{\om}\left(2(n+2)(3G(u)+c^2)+(3n\D u-(n-4)c)(\D
u-c)\right)g(u)dx
\\ & & = 4(n+2)\int_{\om}(\D u-c)\left\{|\na^3 u|^2 -\fr 3{n+2}|\na(\D
u)|^2\right\}dx\\ \no & & \ \ \ +(n+2)\int_{\om}|\na u|^2\D(g(u))dx.
 \en
Here, $\na^3 u=\na(\na^2u)$ is the covariant derivative of the
Hessian $\na^2 u$ of $u$.
\end{theorem}

{\it Another proof of Theorem 3.}  We have from  Rellich identity
that
 \be \int_{\om} udx= -\fr{nc^2|\om|}{n+4}. \en
 Observe that
\be \int_{\om}\D udx=0,\ \int_{\om}(\D u)^2dx=\int_{\om}
u\D^2udx=-\int_{\om} u dx. \en Taking $g(u)=1, \ G(u)=u$, the left
hand side of (1.20) then becomes \be\label{th3.0} \no &
&\int_{\om}\left(2(n+2)(3u+c^2)+(3n\D u-(n-4)c)(\D
u-c)\right)dx\\
&=& \int_{\om}(3(n+4)u+3nc^2)dx=0. \en
 Since $\D^2
u=-1 \ {\rm in}\ \om, \ \D u|_{\pa\om}=c$, we know that $\D u-c>0$
in the interior of $\om$. Therefore, we have from (1.20) and (1.23)
that \be\sum_{i,j,k}u_{ijk}^2-\fr 1{n+2}|\na(\D u)|^2=
|\na^3u|^2-\fr 1{n+2}|\na(\D u)|^2 = 0\en in the interior of $\om$,
and so on $\overline{\om}$ by continuity. Theorem 3 follows as
above. \vskip0.3cm We shall also prove the following symmetry result
using integral identities.
\begin{theorem} Let $\om$ be a bounded domain in $\R^n, n\geq 2$,
with $C^2$ boundary and if the following overdetermined problem has
a solution in $C^2(\overline{\om})$: \be
& & \D u=-1 \ \ \ {\rm in}\   \om, \ \ \ \ u=0\ {\rm on}\ \pa\om, \\
& &   \fr{\pa u}{\pa \nu}= c|x| \ \ \ \ \ \ \  {\ \rm on\ } \pa\om,
\en where  $c$ is a constant. Then $\om$ is a ball centered at the
origin, $c=-\fr 1n$ and \be u(x) = -\fr 1{2n}\left(|x|^2-R^2\right),
\en where $R$ is the radius of $\om$.
\end{theorem}
 When $\om$ contains the origin strictly in its
interior, Theorem 5 has been proven by Tewodros \cite{t} using the
maximum principle.

\section{Proof of the results}
\setcounter{equation}{0}
 In this section, we prove Theorems 2, 4 and 5. Firstly we
make some convention about notation to be used. Let $x=(x_1,\cdots,
x_n)$ and  $\langle, \rangle $ be the position vector and the
standard inner product of $\R^n$, respectively. We shall use $u_i$,
$u_{ij}$, $u_{ijk}$, $u_{ijkl}$ and $u_{ijklm}$ to denote,
respectively, \be\no\fr{\pa u}{\pa x_i},\ \fr{\pa^2 u}{\pa x_i\pa
x_j},\ \fr{\pa^3 u}{\pa x_i\pa x_j\pa x_k},\ \fr{\pa^4 u}{\pa x_i\pa
x_j\pa x_k\pa x_l}\ \ {\rm and}\ \ \fr{\pa^5 u}{\pa x_i\pa x_j\pa
x_k\pa x_l\pa x_m}. \en
\begin{lemma}\label{le1}
Let $u$ satisfy (1.3)-(1.5) and $\eta$ be the solution of the
Dirichlet problem \be\label{le1.1}
 \left\{\begin{array}{l} \D\eta = \lan\na(\D u),\na(\D^2 u)\ra\ \hskip0.03cm \ {\rm in\ \ } \om,\\
 \ \ \hskip0.04cm \eta = 0 \ \ \ \ \ {\rm \ \ \ \ \ \ \ \  \ \ \ \  \ \ \ \ \ on\ \ } \pa\om.
\end{array}
\right. \en The function \be\no F: &=& \fr 12\sum_{i,j,k} u_{ijk}^2
-\fr 12\sum_{i,j} (\D u)_{ij}u_{ij}+\fr 14\lan\na u, \na(\D^2u)\ra
 +\fr{n-8}{4(n+4)}|\na(\D u)|^2\\ & & +\fr 3{(n+4)(n+2)}(\D^2u \D u +u)-\fr{3n(n-2)}{4(n+4)(n+2)}\eta \en
 assumes its maximum value
on $\pa \om$. \end{lemma}

{\it Proof of  Lemma \ref{le1}.} We need only to show that $\D F\geq
0$ in $\om$. A straight forward calculation gives
\be\label{le1.1}\no \D F &=&\sum_{i,j,k,l}
u_{ijkl}^2+\sum_{i,j,k}(\D u)_{ijk} u_{ijk}\\ \no & & -\fr
12\left(\sum_{i,j}\left((\D^2u)_{ij}u_{ij}+(\D
u)_{ij}^2\right)+2\sum_{i,j,k}(\D u)_{ijk} u_{ijk}\right)
\\ \no & & + \fr 14\left( 2\sum_{i,j}(\D^2 u)_{ij}u_{ij}+\lan\na(\D^2 u),\na (\D
u)\ra\right)\\ \no & & +\fr{n-8}{2(n+4)} \left(\sum_{i,j}(\D
u)_{ij}^2+\lan\na (\D^2 u),\na(\D u)\ra\right)
\\ \no & & +\fr 3{(n+4)(n+2)}\left((\D^2 u)^2+ 2\lan\na (\D^2 u),\na(\D
u)\ra\right)\\ \no & & -\fr{3n(n-2)}{4(n+4)(n+2)}\lan\na(\D u),\na(\D^2 u)\ra\\
 &=&\sum_{i,j,k,l} u_{ijkl}^2-\fr 6{n+4}\sum_{i,j}(\D
u)_{ij}^2+\fr{3}{(n+4)(n+2)}(\D^2 u)^2. \en To see that the right
hand side of (\ref{le1.1}) is nonnegative, it suffices to note that
\be\no & & \sum_{i,j,k,l}\left\{u_{ijkl}-\fr 1{n+4}\left((\D
u)_{ij}\de_{kl}+ (\D
u)_{il}\de_{jk}+(\D u)_{ik}\de_{jl}+(\D u)_{jk}\de_{il}\right.\right.\\
\no\label{le1.0} & &\left.\left.\ \ +(\D u)_{jl}\de_{ik}+(\D
u)_{kl}\de_{ij}\right) +\fr {\D^2 u}{(n+4)(n+2)}
(\de_{ij}\de_{kl}+\de_{il}\de_{jk}+\de_{ik}\de_{jl})\right\}^2\\
& & \ \ \ \ \ \ \ \ =\sum_{i,j,k,l} u_{ijkl}^2-\fr
6{n+4}\sum_{i,j}(\D u)_{ij}^2 +\fr{3}{(n+4)(n+2)}(\D^2 u)^2. \en
This completes the proof of Lemma 2.1. \vskip0.3cm {\bf Lemma 2.}
{\it Let $u$ be a solution of (1.3)-(1.5). The following identities
hold: \be\label{le2.0}  \int_{\om} udx=\fr {nc^2|\om|}{n+6}, \en
\be\label{le2.1}\int_{\om} Fdx =
\fr{3(n+2)nc^2|\om|}{2(n+4)(n+6)}-\fr{3n(n-2)}{8(n+4)(n+2)}\int_{\om}(\D^2u)^2\gamma
dx.\en } {\it Proof of Lemma 2.}  It follows from (1.4)  that \be
\na^2 u=0 \ \ {\rm on \ \ }\pa\om. \en Here $\na^2 u $ denotes the
Hessian of $u$ and is given by \be\label{th2peq2} \na^2u(\alpha,
\beta)=\langle \na_{\alpha}\na u, \beta\rangle\en for all $\alpha,
\beta\in \mathfrak{X}(\om).$ From (1.3), we have \be \label{ple2.1}
\D^3\lan x, \na u\ra = 6\D^3 u+ \lan x, \na(\D^3u)\ra =-6. \en
Multiplying (\ref{ple2.1}) by $u$ and integrating on $\om$, one gets
from (1.3)-(1.5), (2.7) and the divergence theorem that
\be\label{ple2.2}\no -6\int_{\om} u dx&=& \int_{\om} u \D^3\lan x,
\na u\ra dx
\\ \no
&=& \int_{\om} \D u\ \D^2\lan x, \na u\ra dx\\ \no &=& -\int_{\om}
\lan \na(\D u), \na(\D\lan x, \na u\ra)\ra dx\\ \no &=& \int_{\om}
\D^2 u\ \D\lan x, \na u\ra dx - \int_{\pa\om}\D\lan x, \na u\ra
\fr{\pa(\D u)}{\pa \nu}ds\\ \no &=& \int_{\om} \D^2 u\ \D\lan x, \na
u\ra dx- c\int_{\pa\om}(2\D u+ \lan x, \na(\D u)\ra)ds\\ \no &=&
\int_{\om}
\D^2 u\ \D\lan x, \na u\ra dx- c\int_{\pa\om}\lan x, \na(\D u)\ra ds\\
\no &=& \int_{\om} \D^2 u\ \D\lan x, \na u\ra dx- c\int_{\pa\om}\lan
x, \nu\ra\fr{\pa(\D u)}{\pa\nu }ds\\ \no &=& \int_{\om} \D^2 u\
\D\lan x, \na u\ra dx- c^2\int_{\pa\om}\lan x, \nu\ra ds
\\ \no &=& -\int_{\om}\lan \na(\D^2 u),
\na\lan x, \na u\ra\ra dx +\int_{\pa\om} \D^2 u \fr{\pa\lan x, \na
u\ra}{\pa\nu}ds -nc^2 |\om|\\ \no &=& \int_{\om}\D^3 u \lan x, \na
u\ra)dx +\int_{\pa\om} \D^2 u(\lan \nu, \na u\ra + \na^2 u(x, \nu))ds- nc^2|\om|\\ \no &=& - \int_{\om}\lan x, \na u\ra dx- nc^2|\om|\\
\no &=& n\int_{\om} u dx-nc^2 |\om|. \en This  proves (\ref{le2.0}).
In order to obtain (\ref{le2.1}), we integrate \be \fr 12\sum_{i,j}
\D(u_{ij}^2)=\sum_{i,j,k} u_{ijk}^2+\sum_{i,j}(\D u)_{ij} u_{ij} \en
on $\om$ and use $u_{ij}|_{\pa\om}=0, \ \forall i, j$, to obtain
\be\label{ple2.3} \sum_{i,j,k}\int_{\om} u_{ijk}^2 dx=
-\sum_{i,j}\int_{\om} (\D u)_{ij} u_{ij}dx. \en Similarly, one gets
by integrating \be \D\lan\na(\D u),\na u\ra= 2\sum_{i,j}(\D u)_{ij}
u_{ij} + |\na(\D u)|^2 +\lan\na(\D^2 u),\na u\ra \en on $\om$ that
\be\label{ple2.4} -\int_{\om} (\D u)_{ij} u_{ij} dx&=&\fr
12\int_{\om}|\na(\D u)|^2dx
 + \fr 12\int_{\om} \lan\na(\D^2 u),\na u\ra dx\\ \no
 &=& - \fr 12\int_{\om} \D u \D^2 u dx- \fr 12\int_{\om} u \D^3 u dx
\\ \no
 &=& \int_{\om}  u dx.
 \en
Integrating $F$ on $\om$ and using (\ref{ple2.3}), (\ref{ple2.4}),
(1.3), (2.5) and the divergence theorem, one has \be\label{ple2.5}
\no \int_{\om} F dx&=&\fr{3(n+2)}{2(n+4)}\int_{\om} u dx
-\fr{3n(n-2)}{4(n+4)(n+2)}\int_{\om}\eta dx\\
&=&\fr{3(n+2)nc^2|\om|}{2(n+4)(n+6)}-\fr{3n(n-2)}{4(n+4)(n+2)}\int_{\om}\eta
dx. \en To finish the proof of (2.6), we need to calculate
$\int_{\om}\eta$. Multiplying the equation \be\no\ \D\eta =
\lan\na(\D u),\na(\D^2 u)\ra \en by $\gamma$ and integrating on
$\om$, we infer \be\no -\int_{\om}\eta dx&=&\int_{\om}\eta\D\gamma dx\\
\no &=& \int_{\om}\gamma\D\eta dx\\ \no
&=&\int_{\om}\gamma\lan\na(\D u),\na(\D^2 u)\ra dx\\ \no &=&
-\int_{\om}\D u\left(\lan\na\gamma,\na(\D^2u)\ra
+\gamma\D^3u\right)dx\\ \no &=&-\int_{\om}\D
u\lan\na\gamma,\na(\D^2u)\ra dx
 +\int_{\om}\gamma\D u dx\\ &=&-\int_{\om}\D
 u\lan\na\gamma,\na(\D^2u)\ra dx
 -\int_{\om} u dx. \en On the other hand, we have \be\no
\int_{\om}\gamma(\D^2u)^2dx&=& \int_{\om}\D u\D(\gamma\D^2u)dx\\ \no
&=& \int_{\om}\D u((\D\gamma)\D^2 u+\gamma\D^3u
+2\lan\na\gamma,\na(\D^2u)\ra)dx\\ \no &=&
 \int_{\om}\D u(-\D^2 u-\gamma
+2\lan\na\gamma,\na(\D^2u)\ra)dx\\ \no &=&
 \int_{\om}u(-\D^3 u-\D\gamma)dx
+2\int_{\om}\D u\lan\na\gamma,\na(\D^2u)\ra dx\\  &=&
 2\int_{\om}u dx
+2\int_{\om}\D u\lan\na\gamma,\na(\D^2u)\ra)dx.\en Combining the
above two equalities, we arrive at \be\label{ple2.6} \int_{\om}\eta
dx = \fr 12\int_{\om}\gamma(\D^2u)^2 dx. \en Substituting (2.17)
into (2.14), we obtain (2.6).

\vskip0.3cm {\it Proof of Theorem 2.} One knows from (1.4) and (1.5)
that \be \sum_{i,j,k}u_{ijk}^2|_{\pa\om}=(|\na(\D
u)|_{\pa\om})^2=\left(\left.\fr{\pa(\D u)}{\pa
\nu}\right|_{\pa\om}\right)^2= c^2. \en Hence, \be F|_{\pa\om}=\fr
{3nc^2}{4(n+4)}, \en which, in turn implies from  Lemma 1 that
\be\label{ple2.5} F\leq \fr{3nc^2}{4(n+4)} \ \ {\rm in }\ \
\overline{\om}. \en When $n=2$, we know from (2.6) that \be\no
\int_{\om} Fdx=\fr{3\cdot 2\cdot c^2|\om|}{4\cdot 6}\en and when
$n\geq 3$, we have from (1.6) and (2.6) that  \be \int_{\om} F
dx\geq \fr{3nc^2|\om|}{4(n+4)}.\en Hence, for $n\geq 2$, $F\equiv
\fr{3nc^2}{4(n+2)}$ in $\overline{\om}$ and so $\D F$ vanishes
identically in $\overline{\om}$. Therefore, each term of the sum on
the left hand side of  (\ref{le1.0}) vanishes. Consequently, we have
 \be\label{pth1.5}\no
u_{ijkl} &=&\fr 1{n+4} \left\{(\D u)_{ij}\de_{kl}+ (\D
u)_{il}\de_{jk}+(\D u)_{ik}\de_{jl}\right. \\ \no & & \left.\ \ \ \
+(\D u)_{jk}\de_{il}+(\D u)_{jl}\de_{ik}+(\D u)_{kl}\de_{ij}\right\}
\\  & & -\fr {\D^2 u}{(n+4)(n+2)} (\de_{ij}\de_{kl}+\de_{il}\de_{jk}+\de_{ik}\de_{jl}), \ \ \forall i, j, k, l.
\en By differentiating the above equality with respect to $x_m$ and
adding, we obtain \be\no\label{pth1.6} \sum_l u_{ijkll}&=&(\D
u)_{ijk}\\  &=&\fr 1{n+2}((\D^2u)_i \de_{jk}+(\D^2u)_j
\de_{ik}+(\D^2u)_k \de_{ij}) \en Differentiating   with respect to
$x_k$ and summing over $k$, one gets \be\label{pth1.} (\D^2
u)_{ij}=-\fr 1n \de_{ij}. \en Thus we have \be\label{pth1.3} \D^2
u(x)=\fr 1{2n}(A-|x-a_0|^2) \en where $A$ is a constant  and $\D^2
u(a_0)=\fr{A}{2n}$. Without loss of generality, we assume that $a_0$
is the origin. Substituting (2.25) into (2.23), we get
\be\label{pth1.7} (\D u)_{ijk}(x)=-\fr 1{n(n+2)}(x_i \de_{jk}+x_j
\de_{ik}+x_k \de_{ij}). \en
Differentiating (\ref{pth1.5}) with
respect to $x_m$ and using (\ref{pth1.7}) and
$$
(\D^2u)_m=-\fr {x_m}n,
$$
we get \be\label{pth1.8} \no
 u_{ijklm}&=&-\fr 1{n(n+2)(n+4)}\left\{x_{i}(\de_{jk}\de_{lm}+ \de_{jl}\de_{km}+\de_{jm}\de_{kl})\right.\\ \no & &
 \ \ \ \ +x_{j}(\de_{ik}\de_{lm} + \de_{il}\de_{km}+\de_{im}\de_{kl})\\ \no & &
 \ \ \ \ +x_{k}(\de_{ji}\de_{lm}+ \de_{jl}\de_{im}+\de_{jm}\de_{il})
 \\ \no & & \ \ \ \ +x_{l}(\de_{jk}\de_{im}+ \de_{ji}\de_{km}+\de_{jm}\de_{ki})
 \\  & & \ \ \ \ \left.+x_{m}(\de_{jk}\de_{li}+ \de_{jl}\de_{ki}+\de_{ji}\de_{kl})\right\}, \ \ \forall i, j, k, l, m.
\en Consider the function $q: \om\ri\R$ given by
$$ q(x)= u(x)+\fr 1{48n(n+2)(n+4)}|x|^6.$$
Using a straightforward calculation and (2.27), we get
$$
q_{ijklm}=0, \ \ \forall i, j, k, l, m.
$$
Thus $q$ is a  polynomial of  $x_1,\cdots, x_n$ of order $4$   and
so \be\label{pth1.9} \D u(x)=-\fr 1{8n(n+2)}|x|^4+p(x).\en Here, $p$
is a quadratic polynomial of $x_1,\cdots, x_n$. Now let us determine
$p$. From (1.4) we know from the divergence theorem that
\be\label{pth1.10} \int_{\om} (\D u) h dx=0 \ \ for\ all\ harmonic\
h\ in\ \om. \en After some calculations by using (\ref{pth1.3}) we
see that
$$ h=\left( x_i\fr{\pa}{\pa x_j}-x_j\fr{\pa}{\pa x_i}\right)\left(x_i(\D
u)_j-x_j(\D u)_i\right)
$$
is harmonic in $\om$. Then integration by parts using $\D
u|_{\pa\om}=0$ results in \be 0&=&\int_{\om}(\D u)\left(
x_i\fr{\pa}{\pa x_j}-x_j\fr{\pa}{\pa x_i}\right)\left(x_i(\D
u)_j-x_j(\D u)_i\right)dx\\ \no &=& -\int_{\om} \left(x_i(\D
u)_j-x_j(\D u)_i\right)^2dx. \en Hence $ x_i(\D u)_j-x_j(\D
u)_i\equiv 0$ in $\om$ and so $\D u$ is a radial function.
Consequently, we have \be\label{pth1.11} \D u(x)=-\fr
1{8n(n+2)}|x|^4 +\kappa_1 |x|^2+\kappa_2,\en where $\kappa_1,
\kappa_2$ are constants. Since $\D u= 0$ on $\pa\om$, $\om$ is a
ball. We note from (\ref{pth1.11}) that \be \label{pth1.12}
\D(x_iu_j-x_ju_i)=0\ \ {\rm in\ } \om \en
 and from (1.4) that
 \be (x_iu_j-x_ju_i)|_{\pa\om}=0.
 \en
 Hence, $x_iu_j-x_ju_i=0$ in $\om$ and so $u$ is a radial function,
 which, combining with (2.31) and the fact that $u$ is a polynomial, gives
 \be
 u(x)= -\fr 1{48n(n+2)(n+4)}|x|^6 +\fr{\kappa_1}{4(n+2)}|x|^4+\fr
 {\kappa_2}{2n}|x|^2+\kappa_3,
 \en
where $\kappa_3$ is a constant. Let us denote by $\rho$ the radius
of $\om$. One deduces  from (1.4) and (1.5) that \be & & -\fr
1{48n(n+2)(n+4)}\rho^6 +\fr{\kappa_1}{4(n+2)}\rho^4+\fr
 {\kappa_2}{2n}\rho^2+\kappa_3=0,  \\ & & -\fr 1{8n(n+2)(n+4)} \rho^4 +\fr
{\kappa_1}{n+2}\rho^2+\fr{\kappa_2}n=0,\\ & & -\fr 1{8n(n+2)}\rho^4
+\kappa_1
\rho^2+\kappa_2=0,\\
 & &
 -\fr 1{2n(n+2)}\rho^3 + 2\kappa_1\rho =c.
\en Solving (2.35)-(2.38), we obtain \be & & \rho
=\left(|c|n(n+2)(n+4)\right)^{\fr 13}, \\
& & \fr {\kappa_1}{4(n+2)}=\left(\fr{c^2}{n(n+2)(n+4)}\right)^{\fr
13}\cdot \fr 1{16},\\ & & \fr{\kappa_2}{2n} =
-\left(c^4n(n+2)(n+4)\right)^{\fr 13}\cdot\fr 1{16},\\ & & \kappa_3
= \fr{c^2n(n+2)(n+4)}{48}. \en Substituting (2.39)-(2.42) into
(2.34), we get (1.8). This completes the proof of Theorem 2.

\vskip0.3cm {\bf Remark.} From (2.16), we have \be\no
\int_{\om}\gamma(\D^2u)^2dx&=& 2\int_{\om} u dx + 2\int_{\om}
u\D\lan\na\gamma, \na(\D^2u)\ra dx\\ &=& 2\int_{\om} u dx +
4\int_{\om} u\left\{\sum_{i,j} \gamma_{ij}(\D^2 u)_{ij}\right\}dx.
\en In the case that $\om$ is a ball with center $a$ and  radius
$R$, $\gamma$ is given by \be \gamma(x)=-\fr{|x-a|^2-R^2}{2n}. \en
Thus \be \gamma_{ij}=-\fr 1n\delta_{ij}, \  \forall i,j, \en which
gives \be \int_{\om} u\left\{\sum_{i,j} \gamma_{ij}(\D^2
u)_{ij}\right\}dx &=& -\fr 1n\int_{\om}u\D^3u =\fr 1n\int_{\om} u
dx. \en We then obtain from (2.43) and (2.5) that \be
\int_{\om}\gamma(\D^2u)^2dx=\left(2+\fr 4n\right)\int_{\om}u dx
=\fr{2(n+2)c^2|\om|}{n+6}. \en That is, (1.6) becomes an equality
when $\om$ is a ball.
 \vskip0.3cm {\it Proof of
Theorem 4.} From (1.18) and (1.19) we know that \be\label{pth4}
|\na^2u|^2 =c^2 \ \ \ {\rm on}\ \ \pa\om. \en
 Multiplying
(1.17) by $|\na^2 u|^2$ and integrating on $\om$, we have
\be\label{pth4.0} -\int_{\om} g(u)|\na^2 u|^2
dx=\int_{\om}(\D^2u)|\na^2u|^2dx.\en Observe that \be\fr 12 \D |\na
u|^2 = |\na^2 u|^2 +\lan\na u, \na(\D u)\ra. \en Using (1.18),
(1.19), (\ref{pth4}) and the divergence theorem, we have
 \be\no\label{pth4.1} -\int_{\om} g(u)|\na^2
u|^2dx&=&-\int_{\om}g(u)\left(\fr 12\D|\na u|^2-\lan\na u, \na(\D
u)\ra\right)dx\\ \no &=&\fr 12\int_{\om}\lan\na(g(u)), \na|\na
u|^2\ra dx+\int_{\om}\lan\na(G(u)),\na(\D u)\ra dx\\ \no &= &-\fr
12\int_{\om}|\na u|^2\D(g(u))dx-\int_{\om}G(u)\D^2 u dx\\ &=& -\fr
12\int_{\om}|\na u|^2\D(g(u))dx+\int_{\om}G(u)g(u)dx,\en

\be\no\label{pth4.2} \int_{\om}(\D^2u)|\na^2u|^2dx
&=&\int_{\pa\om}|\na^2u|^2\fr{\pa(\D u)}{\pa\nu}ds
-\int_{\om}\lan\na(\D u), \na|\na^2u|^2\ra dx\\ \no &=&
c^2\int_{\om}\D^2u dx + \int_{\om}(\D
u)\D|\na^2u|^2dx-\int_{\pa\om}(\D u)\fr{\pa(|\na^2u|^2)}{\pa\nu}ds\\
&=& -c^2\int_{\om} g(u)dx+\int_{\om}(\D u-c)\D|\na^2u|^2 dx. \en We
have
\be\label{pth4.3}\no\D|\na^2u|^2&=&2\sum_{i,j,k}u_{ijk}^2+2\sum_{i,j}u_{ij}(\D
u)_{ij}\\ &=&2|\na^3u|^2+2\sum_{i,j}u_{ij}(\D u)_{ij},\en
 \be\no
\D\lan\na u,\na(\D u)\ra &=& 2\sum_{i,j}u_{ij}(\D u)_{ij}
 +|\na(\D u)|^2+\lan\na u, \na(\D^2u)\ra
 \\ &=& 2\sum_{i,j}u_{ij}(\D u)_{ij}
 +|\na(\D u)|^2-\lan\na u, \na(g(u))\ra,
 \en
Combining (2.52)-(2.54), we get \be\no\label{pth4.2} &
&\int_{\om}(\D^2u)|\na^2u|^2dx\\ \no &=& -c^2\int_{\om}
g(u)dx+\int_{\om}(\D u-c)\left(2|\na^3u|^2 + \D\lan\na u,\na(\D
u)\ra-|\na(\D u)|^2\right.\\ \no & & \hskip4.5cm \left.+\lan\na u,
\na(g(u))\ra\right)dx\\ \no
 &=& -c^2\int_{\om} g(u)dx+2\int_{\om}(\D u-c)\left(|\na^3u|^2-\fr
3{n+2}|\na(\D u)|^2\right)dx\\ \no & & + \int_{\om}(\D
u-c)\left(\D\lan\na u,\na(\D u)\ra+\fr{4-n}{n+2}|\na(\D
u)|^2+\lan\na u, \na(g(u))\right)dx. \en Observing
 \be\no \int_{\om}|\na(\D
u)|^2dx&=&\int_{\pa\om}\D u\fr{\pa(\D u)}{\pa\nu}ds-\int_{\om}\D
u\D^2u dx\\ \no &=& c\int_{\om}\D^2u dx-\int_{\om}\D u\D^2u dx
\\ &=& \int_{\om}(\D u-c)g(u)dx,
\en \be\no \int_{\om}\D u |\na(\D u)|^2dx&=&\fr
12\int_{\om}\lan\na(\D u),\na(\D u)^2\ra dx\\ \no &=& \fr
12\left\{\int_{\pa\om}(\D
u)^2\fr{\pa(\D u)}{\pa\nu}ds-\int_{\om}(\D u)^2\D^2 u dx\right\}\\
&=& \fr 12\int_{\om}\left((\D u)^2-c^2\right) g(u)dx, \en \be\no
\int_{\om}(\D u-c)\D\lan\na u,\na(\D u)\ra dx&=&\int_{\om}\D(\D
u-c)\lan\na u,\na(\D u)\ra dx\\ \no &=& \int_{\om}\D^2u\lan\na
u,\na(\D u)\ra dx\\ \no &=& -\int_{\om} G(u)g(u)dx, \en \be\no
\int_{\om}(\D u-c)\lan\na u, \na(g(u))\ra dx&=&
-\int_{\om}g(u)\left(\lan\na(\D u), \na u\ra +(\D u-c)\D u\right)dx\\
&=& -\int_{\om}\left((\D u-c)\D u+ G(u)\right)g(u)dx, \en one
arrives at
\be\no\int_{\om}(\D^2u)|\na^2u|^2 dx &=&
2\int_{\om}(\D u-c)\left(|\na^3u|^2-\fr 3{n+2}|\na(\D u)|^2\right)dx\\
 & &
-\int_{\om}\left(G(u)+c^2\right)g(u)dx\\
\no & & +\fr{4-n}{2(n+2)}\int_{\om}(\D u-c)^2 g(u)dx\\ \no & &
-\int_{\om}\left((\D u-c)\D u+ G(u)\right)g(u)dx. \en Combining
(2.49), (2.51) with (2.58), we obtain (1.13). This completes the
proof of Theorem 4.

\vskip0.3cm {\it Proof of Theorem 5.} We have \be\int_{\om}|\na
u|^2dx=-\int_{\om} u\D u dx=\int_{\om}u dx.\en Multiplying $\D u=-1$
by $|\na u|^2$ and integrating on $\om$, we get \be\no -\int_{\om} u
dx&=& -\int_{\om}|\na u|^2 dx\\ \no &=&\int_{\om}|\na u|^2\D u dx\\
\no &=& -\int_{\om}\lan\na|\na u|^2, \na u\ra dx + \int_{\pa\om}|\na
u|^2\fr{\pa u}{\pa\nu}ds\\ \no &=& \int_{\om} u\D |\na
u|^2dx+c^3\int_{\pa\om}|x|^3ds\\ &=& 2\int_{\om} u|\na^2 u|^2dx
+c^3\int_{\pa\om}|x|^3ds. \en Multiplying $\D u=-1$ by $\lan x, \na
u\ra$ and integrating on $\om$, one has \be\no n\int_{\om} u dx&=&
-\int_{\om}\lan x, \na u\ra dx\\ \no &=& \int_{\om} \D u\lan x, \na
u\ra dx\\ \no &=& -\int_{\om}\lan\na u, \na\lan x, \na
u\ra\ra dx +\int_{\pa\om}\lan x, \na u\ra\fr{\pa u}{\pa\nu}ds\\
\no &=& \int_{\om} u\D\lan x, \na u\ra dx+\int_{\pa\om}\lan x,
\nu\ra\left(\fr{\pa u}{\pa\nu}\right)^2ds
\\ \no &=& -2\int_{\om} u dx+\int_{\pa\om} c^2\lan x, \nu\ra |x|^2ds.
\\ \no &=& -2\int_{\om} u dx+
\fr{c^2}4\int_{\pa\om}\fr{\pa|x|^4}{\pa\nu}ds
\\ \no &=& -2\int_{\om} u dx+ \fr{c^2}4\int_{\om}\D |x|^4dx
\\ &=& -2\int_{\om} u dx+ c^2(n+2)\int_{\om} |x|^2dx.\en
Thus, we have \be \int_{\om} u dx= c^2\int_{\om} |x|^2dx. \en On the
other hand, multiplying $\D u=-1$ by $|x|^2$ and integrating on
$\om$, one arrives at \be\no -\int_{\om} |x|^2 dx&=& \int_{\om} |x|^2 \D u dx\\
\no &=& -\int_{\om} \lan\na|x|^2, \na u\ra dx+\int_{\pa\om}
|x|^2\fr{\pa u}{\pa\nu}ds\\ \no &=& \int_{\om} u\D |x|^2dx
+c\int_{\pa\om}|x|^3ds\\ &=& 2n\int_{\om} u
dx+c\int_{\pa\om}|x|^3ds. \en Substituting (2.63) into (2.62), we
infer \be (1+2nc^2)\int_{\om}| u dx+ c^3\int_{\pa\om} |x|^3 ds=0.
\en Also, we have \be\no |\om|&=&-\int_{\om} \D u dx=
-\int_{\pa\om}\fr{\pa u}{\pa\nu}ds= (-c)\int_{\pa\om}|x|ds\\ &=&
|c|\int_{\om}|x|dx\geq|c|\int_{\pa\om} \lan x, \nu\ra ds=|c| n|\om|,
\en which gives \be |c|\leq\fr 1n. \en Thus, we have from (2.64) and
(2.65) that \be \left(1+\fr 2n\right)\int_{\om} u dx+
c^3\int_{\pa\om} |x|^3 ds\geq 0. \en It follows from (2.60) and
(2.67) that \be 0\geq\int_{\om} u\left(|\na^2u|^2-\fr
1n\right)dx=\int_{\om} u\left(|\na^2u|^2-\fr{(\D u)^2}n\right)dx.
\en Since $\D u=-1$ in $\om$, $u|_{\pa\om}=0$, $u>0$ in the interior
of $\om$. The Schwarz inequality implies that \be |\na^2u|^2-\fr{(\D
u)^2}n\geq 0.\en Therefore, we conclude from (2.68)  that \be
|\na^2u|^2-\fr{(\D u)^2}n=0\en and that the inequality (2.65) is
actually an  equality. Consequently, we have \be c=-\fr 1n,\en \be
 \ \
u_{ij}=-\fr 1n \delta_{ij}, \ \forall i,j \en and \be x=|x|\nu {\ \
\ \rm on}\ \ \pa\om.\en Consider the function $\beta:\pa\om\ri\R$
given by $\beta(x)=|x|^2$. For any $w\in\mathfrak{X}(\pa\om)$, it
follows from (2.74) that \be w\beta=2\lan x, w\ra=2\lan |x|\nu, w\ra
=0, \en which shows that $\beta$ is a constant. Hence, $\pa\om$ is a
sphere centered at the origin and so $\om$ is a ball centered at the
origin. One then knows from (2.72) that \be u=-\fr
1{2n}(|x|^2-R^2),\en where $R$ is the radius of the ball. This
completes the proof of Theorem 5.

\vskip0.3cm
 \noindent Changyu Xia (xiachangyu666@163.com)


\begin{thebibliography}{23}
\bibitem{a} G. Alessandrini, A symmetry theorem for condensers,
Math. Methods Appl. Sci. 15 (1992), 315-320.
\bibitem{bnst1} B. Brandolini, C.
Nitsch, P. Salani, and C. Trombetti, Serrin type overdetermined
problems: an alternative proof, Arch. Rat. Mech. Anal. 190 (2008),
267-280.
\bibitem{bnst2} B. Brandolini, C. Nitsch, P. Salani, and C. Trombetti, On
the stability of the Serrin problem, J. Diff. Equations 245 (2008),
1566-1583.

\bibitem{b} A. Bennett, Symmetry in an overdetermined fourth order
elliptic boundary value problem, SIAM J. Math. Anal. 17 (1986),
1354-1358.
\bibitem{bnst}B. Brandolini, C. Nitsch, P. Salani, C. Trombetti,
Serrin-Type overdetermined problems: an alternative proof, Arch.
Rat. Mech. Anal. 190 (2008), 267-280.
\bibitem{ch}M. Choulli and A. Henrot, Use
of the domain derivative to prove symmetry results in partial
differential equations, Math. Nachr. 192 (1998), 91?03.

\bibitem{fk}A. Farina and B. Kawohl,
Remarks on an overdetermined boundary value problem, Calc. Var.
Partial Differential Equations 31 (2008), no. 3, 351?57.

\bibitem{gl} N. Garofalo, J. Lewis, A symmetry result related to
some overdetermined boundary value problems, Amer. J. Math. 111
(1989), 9-33.
\bibitem{p}L. E. Payne, Some remarks on maximum principles, J.
Analyse Math. 30 (1976), 421-433.
\bibitem{ps}L. E. Payne and P. W. Schaefer, Duality theorems in some overdetermined problems, Math. Methods in the
Appl. Sciences, 11 (1989), 805-819.
\bibitem{pr}G. A. Philippin, L. Ragoub, On some second order and fourth order elliptic overdetermined
problems, Z. Angew. Math. Phys. 46 (1995) 188-197.
\bibitem{t} W. C. Troy, Symmetry properties in systems of semilinear
elliptic equations, J. Differential Equations 42 (1981), 400-413.
\bibitem{s} J. Serrin,  A symmetry problem in potential theory, Arch.
Rat. Mech. Anal., 43 (1971), 304-18.
\bibitem{t}A. Tewodros, Two symmetry problems in potential theory, EJDE, 43
(2001),  1-5.
\bibitem{w} H. Weinberger,
Remark on the precedingpaper of Serrin, Arch. Rat. Mech. Anal., 43
(1971), 319-20.

\end{thebibliography}
\end{document}